\documentclass[a4paper]{amsart}
\usepackage[margin=2.5cm]{geometry}
\usepackage{float}
\usepackage{xcolor}
\definecolor{darkgreen}{rgb}{0,0.45,0}
\usepackage[colorlinks,urlcolor=blue,citecolor=darkgreen,linkcolor=darkgreen]{hyperref}

\newtheorem*{prop*}{Proposition}

\title[Subset Takeaway and the Boolean lattice]
{Counterexamples to conjectures about Subset Takeaway\break
and counting linear extensions of a Boolean lattice}

\author{Andries E. Brouwer}
\address{}
\email{aeb@cwi.nl}

\author{J. Daniel Christensen}
\address{Department of Mathematics,
         University of Western Ontario,
         London, Ontario,
         Canada}
\email{jdc@uwo.ca}

\date{May 3, 2017}

\begin{document}



\begin{abstract}
We develop an algorithm for efficiently computing recursively defined
functions on posets.  We illustrate this algorithm by disproving
conjectures about the game Subset Takeaway (Chomp on a hypercube) and
computing the number of linear extensions of the lattice of a 7-cube
and related lattices.
\end{abstract}

\maketitle

\section{Introduction}

In this note we solve two problems on small Boolean lattices.
The first is that of finding the optimal strategy for the game
of Chomp (or Subset Takeaway). The second is that of counting the
number of linear extensions. In both cases the method is
recursive descent on normalized posets.

\medskip
Chomp is a game played on a partially ordered set $P$ with smallest element 0.
A move consists of choosing a nonzero element $x$ of $P$ and removing $x$ and all
larger elements from $P$. Whoever is unable to move loses.
This game has been studied for several families of partially ordered sets.
For some history and discussion, see the web page~\cite{B}.
The name Chomp is due to Martin Gardner~\cite{G}.

If $P$ has a largest element 1, different from 0, then a trivial
strategy-stealing argument shows that the first player wins.
(If choosing 1 does not win, it is because the opponent has the
devastating reply $a$. But in that case the first player wins
by starting with $a$.)
However, this argument is non-constructive: the winning move is unknown.

Subset Takeaway is the special case of Chomp, where $P$ is chosen
to be the Boolean lattice $B_n$ of all subsets of an $n$-set.
It was studied by Gale and Neyman~\cite{GN}.
Since there is a largest element, the first player wins. But how?
Gale and Neyman asked whether taking the top element, the unique $n$-set,
is the winning move, and this is true for $n \le 5$~(\cite{GN})
and for $n = 6$~(\cite{CT}).
More generally, let $P = P_{n,k}$ be the partially ordered set
of all subsets of size at most $k$ from an $n$-set.
Gale and Neyman asked whether the first player loses in $P_{n,k}$
if and only if $(k+1) \mid n$, and prove this for $k \le 2$.

However, we show that these longstanding conjectures are false. Indeed, using computer search,
one finds fairly easily that the first player loses in $P_{7,3}$
and with more effort that the first player wins in $P_{7,6}$
(by choosing any 4-set).
In particular, the top move in $P_{7,7}$ is not a winning move, but a 4-set is.

In Section~\ref{se:sst}, we extend the above results by giving computations
of Grundy values of $P_{n,k}$ for some values of $n$ and $k$, showing that a
natural pattern fails to continue.
We also list some positions with large Grundy values,
and describe how a reduction technique using involutions can extend
our results to some larger positions.
In Section~\ref{se:linext}, we apply our methods to the problem
of counting the number of linear extensions of the poset $B_n$.
The asymptotic behaviour of this sequence of numbers is of interest,
and bounds were given in \cite{BT,SK}.
Our methods of computation are much more efficient than previous
methods, allowing us to easily reproduce known results for $n \leq 6$
and to compute the answer for $n = 7$ for the first time.
In addition, we derive a formula for $e(P_{n,2})$, and give computations of $e(P_{n,3})$.
Finally, in Section~\ref{se:computation}, we discuss the size of
the search space, describe the techniques we used to efficiently
search it, and give some sample running times and space requirements
that could be useful for future work.

We thank Alejandro Morales for bringing to our attention the problem
of counting linear extensions.

\section{Subset Takeaway}\label{se:sst}
\subsection{Win/loss computation}
We computed that the only winning move in $P_{7,7}$ is to choose a set of size~4.
As we will explain below, the search tree is much too large to be analyzed
by hand, but we can give the first two layers.
A move of size 1 leads to $P_{6,6}$, in which the winning move is known to
be the set of size 6.
Symmetrically, the response to a set of size 6 is the complement.
Similarly, the winning response to a move of size 2 or 5 is the complement;
this follows from the symmetry exchanging the elements of the move of size 2,
using the technique described in Subsection~\ref{ss:involutions}.
However, the winning response to a move of size 3 is not the complement.
For example, our programs show that a winning response to $\{0,1,2\}$
is $\{0,1,3,4,5\}$.
The remaining case is a move of size 4, which our programs show is a winning move.

\subsection{Grundy values}
Any impartial game, such as Subset Takeaway, is equivalent in the sense of
combinatorial game theory to the game of Nim~(cf.~\cite{JHC}).
Thus to each game position $P$ one can associate a nonnegative integer $g$
such that $P$ is equivalent to a Nim heap of size $g$.
This $g$ is called the \emph{Grundy value} of $P$, and
is the smallest nonnegative integer that is not the Grundy value
of one of the positions that can be reached by one move.
The position is lost for the first player if and only if
the Grundy value is 0.
We give the Grundy values of $P_{n,k}$ for $0 \le k \le n$ in small cases.

\begin{table}[H]
\begin{tabular}{c|ccccccccc}
$n\backslash k$ & 0 & 1 & 2 & 3 & 4 & 5 & 6 & 7 & 8 \\
\hline
 0 &  0 \\
 1 &  0 & 1 \\
 2 &  0 & 0 & 2 \\
 3 &  0 & 1 & 0 & 3 \\
 4 &  0 & 0 & 1 & 0     & \bf 1 \\
 5 &  0 & 1 & 2 & 1     &     0 & \bf 2 \\
 6 &  0 & 0 & 0 & 2     & \bf 2 &     0 & \bf 3 \\
 7 &  0 & 1 & 1 & \bf 0 & \bf 3 &     1 & \bf 5 & \bf 6 \\
 8 &  0 & 0 & 2 
\end{tabular}
\end{table}
\noindent
In~\cite{KY}, the Grundy values of complete multipartite graphs are determined.
In particular, they show that the Grundy value of $P_{n,2}$ is $n$ modulo $3$,
agreeing with the $k=2$ column of the table.
Note that, without the results in bold, one might conjecture that the
Grundy value of $P_{n,k}$ is $n$ modulo $k+1$.
We believe that most of the results for $k > 2$ are new, except for the
zeros when $k = n - 1$ and $n \leq 6$.
Other values are unknown.

We have seen already that $P_{n,n}$ is a first player win for all $n > 0$.
We shall see below that $P_{8,3}$, $P_{8,5}$, and $P_{9,3}$
are first player wins (and hence have a positive Grundy value).

\subsection{Large Grundy values}
Some positions in these games have rather large Grundy values $g$.
In the table below we give $g$ and the maximal sets of a simplicial complex
with this Grundy value, for various positions with $n = 7$.
We write $ab$ instead of $\{a,b\}$, etc.

\begin{table}[H]
\begin{tabular}{c|l}
$g$ & antichain \\
\hline
9  & 01, 02, 03, 04, 05, 06, 12, 13, 14, 15, 23, 24, 35, 46 \\[3pt]
23 & 012, 013, 014, 015, 023, 024, 025, 026, 034, 045, \\
   & 056, 123, 126, 135, 136, 145, 236, 245, 346, 456 \\[3pt]
37 & 0124, 0134, 0234, 1234, 0125, 0235, 1235, 0145, 0245, 1245, \\
   & 0345, 2345, 0126, 0136, 0236, 1236, 0146, 0246, 1246, 0156, 1356 \\[3pt]
44 & 012345, 01236, 01246, 01346, 2346, 01256, 1356, 2356, 0456
\end{tabular}
\end{table}

\subsection{Involutions}\label{ss:involutions}
Suppose $P$ is a poset and $\phi : P \to P$ an automorphism of order 2
such that $P_0 = \{ x \in P \mid \phi(x) = x \}$ is downwards closed.
Then the games on $P$ and $P_0$ have the same outcome.
(Indeed, if $P_0$ is a second player win, then the second player can answer
each move $x$ outside $P_0$ with $\phi(x)$, and each move inside $P_0$
with the winning strategy for $P_0$. If $P_0$ is a 1st player win,
with winning move $x_0$, then the 1st player can play $x_0$ and then
apply the previous sentence.)
More generally, $P$ and $P_0$ have the same Grundy value (cf.~\cite{KY}).

It is easy to see that $P_{8,3}$ is a first player win: just pick
a singleton, reducing to $P_{7,3}$ which is a first player loss.
This involution argument allows one to see that $P_{8,5}$
is also a first player win: pick a 2-set, say $\{6,7\}$, and apply
the transposition $(67)$ to reduce to $P_{6,5}$ which is
a first player loss. More generally, we see that if $P_{n,k}$
is a first player loss, then $P_{n+1,k}$ is a first player win
if $k > 0$, and $P_{n+2,k}$ is a first player win if $k > 1$.
So, $P_{9,3}$ is also a first player win.

\section{Counting linear extensions of a Boolean lattice}\label{se:linext}
The machinery we developed for studying Subset Takeaway on an $n$-set
can be used more generally to recursively compute a function whose
value on a poset is defined in terms of its values on smaller posets.
Another example of such a function is $e(P)$, the number of linear
extensions of the poset $P$.
A \emph{linear extension} of a poset $P$ is a total order on $P$ which is
compatible with the given partial order.
We can compute the number of linear extensions of $P$ in the following way.
The maximum element of a total order must be one of the maximal elements
of $P$, so $e(P)$ is given by
\[
e(P) = \begin{cases}
         1,                                          & \mbox{if $P = \emptyset$,}\\[5pt]
         \displaystyle\sum_m \, e(P\setminus \{m\}), & \mbox{otherwise,}
       \end{cases}
\]
where the sum is over the maximal elements $m$ of $P$.
The values of $e(P)$ when $P$ is the Boolean lattice $B_n$
of subsets of an $n$-set were known for $n \le 6$
(see OEIS A046873~\cite{OEIS}).
We also give the value for $n = 7$ (roughly $6.3 \cdot 10^{137}$).

\begin{table}[H]
\begin{tabular}{cl}
$n$ & $e(B_n)$ \\
\hline
0 & 1 \\
1 & 1 \\
2 & 2 \\
3 & 48 \\
4 & 1680384 \\
5 & 14807804035657359360 \\
6 & 141377911697227887117195970316200795630205476957716480 \\
7 & 630470261306055898099742878692134361829979979674711225065761605059425- \\
  & 237453564989302659882866111738567871048772795838071474370002961694720
\end{tabular}
\end{table}
\vspace*{-10pt}\hspace*{-2.2cm}
\textcolor{white}
{\tiny 630470261306055898099742878692134361829979979674711225065761605059425237453564989302659882866111738567871048772795838071474370002961694720}

\noindent
The value 70016 of $e(B_4)/4!$ was given by Weinrich~\cite{W}.
The value of $e(B_5)$ was given by Morton~\cite{M},
and $e(B_6)$ was found by Wienand~\cite{Wi}.
Wienand reports that the computation of $e(B_6)$ took less than 16 hours
on a compute server. For us that same computation took less than
1 second, while the computation of $e(B_7)$ took less than 2.5 hours.
Asymptotic bounds on $e(B_n)$ were given in \cite{BT,SK}.
In \cite{BW} it was shown that counting linear extensions is \#P complete.


As a further test of our approach, one can derive the following explicit formula
for $e(P_{n,2})$.

\begin{prop*}
For $n \geq 0$,
\[
e(P_{n,2}) = \frac{n! \, \left(\binom{n}{2}+n\right)!}{\prod\limits_{i=1}^n in - \binom{i}{2}} .
\]
\end{prop*}

This formula gives 183516891399909333860213587968000000
as the value of $e(P_{8,2})$, which agrees with the results
of our programs.

\begin{proof}
The equation holds for $n=0,1$. Let $n \ge 2$.
The factor of $n!$ in the numerator counts the orderings of the $n$ singletons.
Assuming that they have been ordered $1 < 2 < \cdots < n$, we count the possible
placements of the $2$-sets.
The $2$-set $\{n-1,\, n\}$ has one possible place.
The $2$-set $\{n-2,\, n\}$ then has two possible places.
This continues until the $2$-set $\{1,\, n\}$, which has $n-1$ possible places.
The next $2$-set $\{n-2,\, n-1\}$ has $n+1$ places, the larger jump being
because it can appear before or after $n$.
The numbers increase by one until $\{1,\, n-1\}$, and then there is another
jump by two for $\{n-3,\, n-2\}$.
The process ends when we reach $\{1,\, 2\}$, which has $\binom{n}{2} + n - 2$ places.
The result is a product of all $k$ between $1$ and $\binom{n}{2} + n - 2$,
except for those that are of the form
\[
  \sum_{r=0}^{i-1} \, n-r = i n - \binom{i}{2}
\]
for some $1 \leq i \leq n-2$.
Thus
\[
e(P_{n,2}) = \frac{n! \, \left(\binom{n}{2}+n - 2\right)!}{\prod\limits_{i=1}^{n-2} i n - \binom{i}{2}}
          = \frac{n! \, \left(\binom{n}{2}+n\right)!}{\prod\limits_{i=1}^n i n - \binom{i}{2}} .
\vspace*{-10pt}
\]
\end{proof}

Finding a formula for $e(P_{n,3})$ appears to be more difficult.
Our computations of this for $n = 3, 4, 5, 6$ and $7$ give
48, 1680384, 37783650956544000, 567722883627880394131962569026437120 and
1146874473\-4\-4\-9\-4\-754078263804379304906839548713527318697219614310400000.

\section{Computation}\label{se:computation}

\subsection{Counting antichains}
How much work is it to analyze $P_{n,n}$?
The full game tree consists of all downwards closed families of sets
on $n$ points, or, equivalently (by taking maximal elements), of all
antichains on $n$ points. The number of antichains on $n$ points
is known as the $n$-th Dedekind number $M(n)$, cf. OEIS A000372~\cite{OEIS}.
A rough lower bound is found by taking arbitrary families of sets
of size $\lfloor n/2 \rfloor$, and we see that
$\log_2 M(n) \ge \binom{n}{\lfloor n/2 \rfloor}$.
For $0 \le n \le 8$ the values of $M(n)$ are
2, 3, 6, 20, 168, 7581, 7828354, 2414682040998, 56130437228687557907788
(cf.~\cite{Wiedemann}).
Here both the empty antichain and the antichain consisting of just
the empty set are counted.

The amount of work can be reduced by using
the action of the symmetric group ${\rm Sym} (n)$ of order $n!$.
This saves at most a factor of $n!$, and for $n=6, 7$, and 8 the value of
$M(n)/n!$ is approximately 10873, 479103580, and 1392123939203560464.
The number of antichains on $n$ points up to isomorphism, or, equivalently,
of unlabeled simplicial complexes with at most $n$ vertices,
is given in OEIS A003182~\cite{OEIS}.
For $0 \le n \le 7$ the values are 2, 3, 5, 10, 30, 210, 16353, and 490013148.
Note that our approximation $M(n)/n!$ gets better as $n$ increases,
and so we expect something like $1.4 \cdot 10^{18}$ antichains up to isomorphism for $n = 8$.

We see that in order to answer questions about $P_{n,n}$ for $n = 7$
it suffices to look at at most 490013148 posets, which is not too much
for modern computers, while a brute force approach for $n = 8$
requires looking at more than $10^{18}$ posets, which is infeasible today.

\subsection{Methods}

Both authors independently wrote and ran computer programs to verify
the claims, and compared the programs on many test cases and against
known positions.

As discussed above, one wants to use the action of ${\rm Sym} (n)$.
This is done by bringing each downwards closed set into a canonical form
and saving already computed values in a hash table.
Given a downwards closed set $A$ on the set $N = \{0,\ldots,n-1\}$,
we consider the invariant that assigns to each $i \in N$ the vector
whose $j$th entry ($j=1,\ldots,n$) is the number of $j$-sets in $A$
containing $i$.  We then apply a permutation to $N$ that brings these invariants
into some lexmin order.  Finally, we consider all permutations in ${\rm Sym} (n)$
that preserve the invariants and apply the one that puts $A$ in a lexmin form.
Putting moves into canonical form is crucial for time and memory
efficiency, and was also the bottleneck in our programs.

\subsection{Sample running times and memory usage}

The programs were written in C and were run on ordinary desktop hardware
using a single core.
Here we give sample computation times and posets stored.
It would certainly be possible to improve these with more work,
but they give an illustration of the resources we needed,
which may help someone who wishes to reproduce or extend our work.

\begin{table}[H]
\begin{tabular}{clcr}
poset & computation          & time                  & posets \\
      &                      & (minutes)             & stored \\
\hline                        
$P_{7,3}$ & win/loss          & \makebox[2em][r]{5}   & 9680539 \\
$P_{7,7}$ & win/loss          & \makebox[2em][r]{230} & 401448670 \\
\hline                                    
$P_{7,3}$ & Grundy            & \makebox[2em][r]{22}  & 15466911 \\
$P_{7,7}$ & Grundy            & \makebox[2em][r]{637} & 490013147 \\
\hline                                    
$P_{7,3}$ & linear extensions & \makebox[2em][r]{5}   & 15466911 \\
$P_{7,7}$ & linear extensions & \makebox[2em][r]{148} & 490013147 \\
\end{tabular}
\end{table}

For all three types of computation, $P_{6,6}$ took less than 1 second
and stored at most 16352 posets in the hash table.
All of the computations required less than 32 GB of virtual memory,
except for the last one, counting the number of linear extensions of $B_7$.
It required 38 GB, but it would not be hard to reduce this
by ordering the search by the number of elements of the posets.

\end{document}